\documentclass[12pt]{amsart}
\usepackage{amssymb,latexsym,comment,url}

\newcommand{\isom}{ \cong }

\newcommand{\Q}{{\mathbb Q}}

\newcommand{\rank}{\operatorname{rank}}
\newcommand{\genus}{\operatorname{genus}}

\newcommand{\Z}{{\mathbb Z}}

\newenvironment{Proof}{\par\noindent{\sc Proof:}}%
                      {\hspace*{\fill}\nobreak$\Box$\par\medskip}
                       {\hspace*{\fill}\nobreak$\Box$\par\medskip}

\newtheorem{Proposition}{Proposition}[section]
\newtheorem{Theorem}[Proposition]{Theorem}
\newtheorem{Lemma}[Proposition]{Lemma}

\theoremstyle{definition}

\newtheorem{Remark}[Proposition]{Remark}

\renewcommand{\baselinestretch}{1.1}

\begin{document}

\title[Quadratic twists of pairs of Elliptic Curves]%
{High rank quadratic twists of pairs of elliptic curves}

\author[M. Alaa]%
{Mohamed~Alaa}
\address{Department of Mathematics, Faculty of Science, Cairo University, Giza, Egypt}
\email{malaa@sci.cu.edu.eg}

\author[M. Sadek]%
{Mohammad~Sadek}
\address{American University in Cairo, Mathematics and Actuarial Science Department, AUC Avenue, New Cairo, Egypt}
\email{mmsadek@aucegypt.edu}

\begin{abstract}
 Given a pair of elliptic curves $E_1$ and $E_2$ over the rational field $\mathbb Q$ whose $j$-invariants are not simultaneously 0 or 1728, Kuwata and Wang proved the existence of infinitely many square-free rationals $d$ such that the $d$-quadratic twists of $E_1$ and $E_2$ are both of positive rank. We construct infinite families of pairs of elliptic curves $E_1$ and $E_2$ over $\Q$ such that for each pair there exist infinitely many square-free rationals $d$ for which the $d$-quadratic twists of $E_1$ and $E_2$ are both of rank at least 2.
\end{abstract}

\maketitle


\section{Introduction}

Goldfeld Conjecture states that the average rank of elliptic curves over $\Q$ in families of quadratic twists is $1/2$. This reflects the strong belief amongst number theorists that quadratic twists with rank at least $2$ of an elliptic curve defined over $\Q$ are seldom. In fact, one may find a great deal of literature investigating the rank frequencies for quadratic twists of elliptic curves.

In \cite{Mestre, Mestre1}, Mestre proved that given an elliptic curve over $\Q$, there are infinitely many quadratic twists whose rank is at least $2$. Furthermore, he introduced infinitely many elliptic curves with infinitely many quadratic twists whose rank is at least $3$. He proved, moreover, that if $E$ is an elliptic curve over $\Q$ whose torsion subgroup is isomorphic to $\Z/2\Z\times\Z/8\Z$, then it has infinitely many quadratic twists with rank at least $4$.

In \cite{Rubin-Silverberg,Rubin-Silverberg2}, quadratic twists of ranks at least $2$ and $3$ were introduced. In addition, infinitely many elliptic curves with infinitely many quadratic twists of rank at least $4$ were constructed. For the latter families, the quadratic twists are parametrized by an elliptic curve of positive rank.

Many of the families of elliptic curves that were constructed in \cite{Rubin-Silverberg,Rubin-Silverberg2} are families of Legendre elliptic curves. Those are elliptic curves described by the equation $y^2=x(x-1)(x-\lambda)$ where $\displaystyle\lambda\in\Q-\{0,1\}$. In \cite{Dujella}, more families of Legendre elliptic curves with infinitely many quadratic twists of rank at least $2$ were displayed.

 In \cite{Kuwata-Wang}, the study of quadratic twists of pairs of elliptic curves over $\Q$ was initiated. Given a pair of elliptic curves $E_1$ and $E_2$ over the rational field $\mathbb Q$ whose $j$-invariants are not simultaneously 0 or 1728, it was proved that there exist infinitely many square-free rational numbers $d$ such that the $d$-quadratic twists of $E_1$ and $E_2$ are both of positive rank. Similar questions were posed in \cite{Coogan-Jimenez}. Examples of infinitely many pairs of elliptic curves $E_1$ and $E_2$ and infinitely many rational numbers $d$ were given such that $\rank E_1^d=\rank E_2^d=0$, where $E_i^d$ is the quadratic twist of $E_i$ by $d$. Other examples for which $\rank E_1^d=0$ whereas $\rank E_2^d>0$ were given. In \cite{Bo}, it was shown that given four elliptic curves $E_i$, $i=1,2,3,4$, defined over a number field $K$, there exists a number field $L$ containing $K$ such that there are infinitely many $d\in L$ such that the quadratic twist of $E_i$ by $d$ is of positive rank over $L$.

 In this note, we consider the following question: Can we find infinitely many pairs of elliptic curves $E_1$ and $E_2$, and infinitely many rational numbers $d$ such that the quadratic twists of $E_1$ and $E_2$ by $d$ are both of rank at least $2$? We answer the question in the affirmative. In order to achieve this, we start with two Legendre elliptic curves $E_1$ and $E_2$ over $\Q$, then we use techniques from \cite{Rubin-Silverberg,Rubin-Silverberg2} to create a function field over which quadratic twists of $E_1$ and $E_2$ by a certain polynomial are of rank at least $2$. Then we show that there are infinitely many rational numbers at which when we specialize the latter function field we obtain the rational field $\Q$.

 An explicit description of our construction is given as follows. We find infinitely many pairs $E_1$ and $E_2$ of Legendre elliptic curves parametrized by the projective line. To this family of pairs, we associate a polynomial $g(t)$ of degree $12$ such that the quadratic twists $E_1^{g(t)}$ and $E_2^{g(t)}$ are of rank at least $2$ over $\Q$ if $t$ is the $x$-coordinate of a rational point on a specific elliptic curve of positive rank. In other words, the quadratic twists are parametrized by an elliptic curve of rank at least $1$.

\section{Quadratic twists of positive rank}
\label{sec:2}
Throughout this note, if $E$ is an elliptic curve defined over $\Q$ by the equation $y^2=f(x)$, we write $E^d:dy^2=f(x)$ for the quadratic twist of $E$ by $d$.

In the following section, we collect the preliminaries that we are going to use throughout the note. The results introduced here can be found in \cite{Rubin-Silverberg,Rubin-Silverberg2}.

 Given an elliptic curve $E$ over $\Q$, the following lemma, \cite[Corollary 2.2]{Rubin-Silverberg}, can be used to construct quadratic twists of $E$ with large rank over extensions of $\Q(t)$ by square roots of rational functions in $\Q(t)^{\times}$.

\begin{Lemma}
\label{lem1}
Suppose $E$ is an elliptic curve over $\Q$, $g_1,\ldots,g_r\in\Q(t)^{\times}$, the fields $\Q(t,\sqrt{g_i})$ are distinct quadratic field extensions of $\Q(t)$, and $\rank(E^{g_i}(\Q(t)))>0$ for every $i$. Then
\[\rank(E^{g_1}(\Q(t,\sqrt{g_1g_2},\ldots,\sqrt{g_1g_r})))\ge r.\]
If in addition $\Q(t,\sqrt{g_1g_2},\ldots,\sqrt{g_1g_r})=\Q(u)$ for some $u$, and $g(u)=g_1(t)$, then $\rank(E^{g(u)}(\Q(u)))\ge r$.
\end{Lemma}

Lemma \ref{lem1} can be reformulated as follows, see Proposition 2.3 in \cite{Rubin-Silverberg2}.

\begin{Proposition}
\label{prop:1}
Let $E$ be an elliptic curve defined over $\Q$ by $y^2=f(x)$. Let $h_1(t)=t$, and $h_2(t),\ldots,h_r(t)\in\Q(t)$ non-constant. Let $k_i$ be the square-free part of $f(h_i(t))/f(t)$, and suppose $k_1(t),\ldots, k_r(t)$ are distinct modulo $\Q(t)^2$. Then
\begin{itemize}
\item[(i)] $\rank E^{f(t)}(\Q(t,\sqrt{k_2(t)},\ldots,\sqrt{k_r(t)}))\ge r$;
\item[(ii)] if $C$ is the curve defined by the equations $s_i^2=k_i(t)$ for $i=1,\ldots,r$, then for all but at most finitely many rational points $(t,s_1,\ldots,s_r)=(\tau,\sigma_1,\ldots,\sigma_r)\in C(\Q)$, one has $\rank E^{f(\tau)}(\Q)\ge r$.
\end{itemize}
\end{Proposition}

In \cite{Rubin-Silverberg2}, the rational function $h_i$ is chosen to be a linear fractional transformation $\displaystyle \frac{\alpha t+\beta}{t+\delta}$ which permutes the roots of $f$. Furthermore, a direct calculation shows that $f(h_i(t))=f(\alpha)(t+\delta)f(t)(t+\delta)^{-4}$.

As a direct consequence of Lemma \ref{lem1}, one obtains a method to construct quadratic twists with positive rank, see Lemma 2.3 and Remark 2.4 in \cite{Rubin-Silverberg}.
\begin{Lemma}
\label{lem2}
Let $g\in \Q(t)$. Suppose $E$ is an elliptic curve over $\Q$ defined by $y^2=f(x)$. One has $\displaystyle\rank(E^{g}(\Q(t)))>0$ if and only if there is $h\in\Q(t)$ such that $E^g\isom E^{f\circ h}$.
\end{Lemma}
In fact, if $\rank E^g(\Q(t))>0$ where $(h,k)$ is a point of infinite order , then $f\circ h=k^2 g$.

One finishes this section with the following remark providing an upper bound on the rank of quadratic twists.
\begin{Remark}
\label{rem1}
Suppose $g(t)\in\Q[t]$ is square-free and non-constant, and let $C$ be the curve $s^2=g(t)$. Then
\[\rank(E^g(\Q(t)))\le \genus(C)=\lfloor (\deg g-1)/2\rfloor.\]
\end{Remark}

\section{Quadratic twists of pairs of elliptic curves}
\label{sec:3}
In this section, given two elliptic curves $E_1$ and $E_2$ over $\Q$, we start finding conditions for a rational function $g\in\Q(t)$ such that the quadratic twists $E_1^g$ and $E_2^g$ are of positive rank.

One knows that if $E^{g_1}$ and $E^{g_2}$ are quadratic twists of an elliptic curve $E$, then they are isomorphic if and only if $g_1/g_2$ is a square. This together with Lemma \ref{lem2} yield the following result.

\begin{Lemma}
\label{Lem:pairspositiverank}
Given two non-isomorphic elliptic curves $E_i:y^2=f_i(x)$, $f_i(x)\in\Z[x]$, $i=1,2$, and $g\in\Q(t)$, the quadratic twists $E_1^{g}$ and $E_2^{g}$ are of positive rank over $\Q(t)$ if and only if there exist non-constant rational functions $h_i,M_i\in\Q(t)$, $i=1,2$, such that $E_i^g\isom E_i^{M_i^2(f_i\circ h_i)}$.
In particular, if there exist rational functions $h_1$, $h_2$, and $M$ such that $(f_2\circ h_2)/(f_1\circ h_1)=M^2$, then $\rank E_i^{g}(\Q(t))>0$, where $g\equiv f_i\circ h_i$ modulo $\Q(t)^2$.
\end{Lemma}

The following theorem is the main tool that we use to find infinitely many square-free rationals $d$ and infinitely many pairs of elliptic curves over $\Q$ whose $d$-quadratic twists are of Mordell-Weil rank at least 2.

\begin{Theorem}
\label{Thm:1new}
Let $E_1$ and $E_2$ be two non-isomorphic elliptic curves over $\Q$ defined by $y^2=f_i(x)$, $i=1,2$, where $\deg f_i=3$. Let $h_i(t)\in\Q(t)$ be such that $k_{i}(t)$ is the square-free part of $f_i(h_{i}(t))/f_i(t)$, $i=1,2$, and that $k_1(t)$ and $k_2(t)$ are distinct modulo $\Q(t)^2$. Then:
\begin{itemize}
\item[(i)] $\rank E_i^{f_1(t)}\left(\Q\left(t,\sqrt{k_1(t)},\sqrt{k_2(t)},\sqrt{f_2(t)/f_1(t)}\right)\right)\ge 2$, $i=1,2$,
\item[(ii)] if $C$ is the genus one curve \[z_1^2=k_1(t),\;z_2^2=k_2(t),\; z_3^2=f_2(t)/f_1(t),\]
then for all but finitely many rational points $(t,z_1,z_2,z_3)=(t_0,s_1,s_2,s_3)\in C(\Q)$, one has $\rank E_i^{f_1(t_0)}(\Q)\ge 2$, $i=1,2$.
\end{itemize}
\end{Theorem}
\begin{Proof}
In view of Proposition \ref{prop:1}, one knows that $\rank E_i^{f_i(t)}\left(t,\sqrt{k_1(t)},\sqrt{k_2(t)}\right)\ge 2$. Writing $f_2=f_1\left(\sqrt{f_2/f_1}\right)^2$, this implies that $\rank E_2^{f_1(t)}\left(\Q\left(t,\sqrt{k_1(t)},\sqrt{k_2(t)},\sqrt{f_1(t)/f_2(t)}\right)\right)\ge 2$, hence (i) follows. Part (ii) follows from part (i) and Proposition \ref{prop:1} (ii).
\end{Proof}

\begin{Remark}
The two independent infinite points on $E_i^{f_1(t_0)}$ are $(t_0, \sqrt{f_i(t_0)/f_1(t_0)})$ and $\left(h_i(t_0), \sqrt{f_i(h_i(t_0))/f_1(t_0)}\right)$, $i=1,2$.
\end{Remark}
\begin{Remark}
\label{rem2}
 If $h_i(t)$ is a linear fractional transformation of the form $\displaystyle \frac{\alpha t+\beta}{t+\delta}$, then this implies that $k_i(t)$ is the square-free part of $f_i(\alpha)(t+\delta)$. In this paper, we look for a pair $h_1$ and $h_2$ of linear fractional transformations for which there is an infinite rational point on the curve $C$ in Theorem \ref{Thm:1new}. Once one finds $t_0\in\Q$ such that the point $(t,z_1,z_2,z_3)=(t_0,s_1,s_2,s_3)\in C(\Q)$ is an infinite rational point, this will imply the existence of infinitely many such rational points. In other words, there are infinitely many $t\in \Q$ such that $E_i^{f_1(t)}$ has rank at least $2$ over $\Q$, $i=1,2$.
\end{Remark}

\section{Auxiliary results}
\label{sec:4}

 We recall that a Legendre curve is an elliptic curve described by a Weierstrass equation of the form $E:y^2=x(x-1)(x-\lambda)$, $\lambda\in\Q-\{0,1\}$. We recall that two Legendre curves $E_i$ are isomorphic if and only if \[\displaystyle \lambda_2\in \{\lambda_1,\;1-\lambda_1,\;1/\lambda_1,\;1/(1-\lambda_1),\;\lambda_1/(\lambda_1-1),\;(\lambda_1-1)/\lambda_1\}.\]

In this section, we show that for any pair of rational numbers $(\lambda_1,\lambda_2)$, there is a polynomial $g$ of degree $6$ such that the rank of the quadratic twist of $E_{i}:y^2=f_i(x)=x(x-1)(x-\lambda_i)$, $i=1,2$, by $g$ is positive.

 In order to find $g(t)\in\Q(t)$ such that the $g(t)$-quadratic twist of $E_{i}:y^2=f_i(x)=x(x-1)(x-\lambda_i)$ is of positive rank, one finds $h_1,h_2,u\in\Q(t)$ such that $f_2\circ h_2=u^2(f_1\circ h_1)$, see Lemma \ref{Lem:pairspositiverank}. Setting $h_1(t)=h_2(t)=t$, we need to find $t\in\Q$ such that $f_2(t)=u^2 f_1(t)$. In other words, one has $\displaystyle t=\frac{-\lambda_2+\lambda_1u^2}{-1+u^2}$. Therefore, Choosing $g$ to be the square-free part of $f_1(t(u))=f_2(t(u))$ mod $\Q(u)^2$ yields the following result.

 \begin{Lemma}
 \label{lem3}
Let $E_{i}:y^2=f_i(x)=x(x-1)(x-\lambda_i)$, $i=1,2$, be two non-isomorphic Legendre curves. Consider the polynomial
$$g(u)=(\lambda_1-\lambda_2)(-1+u^2)(1-\lambda_2+(-1+\lambda_1)u^2)(-\lambda_2+\lambda_1u^2).$$
Then the Mordell-Weil rank $r$ of the quadratic twists $E_{1}^{g},E_{2}^{g}$ over $\Q(u)$ satisfies $1\le r\le 2$.
\end{Lemma}
\begin{Proof}
It is direct calculation to check that the points $\displaystyle P_{1}=\left(\frac{-\lambda_2+\lambda_1u^2}{-1+u^2},\frac{1}{(-1+u^2)^2}\right)$ lies in $E_{1}^{g(u)}(\Q(u))$ and $\displaystyle P_{2}=\left(\frac{-\lambda_2+\lambda_1u^2}{-1+u^2},\frac{u}{(-1+u^2)^2}\right)$ is a point in $E_{2}^{g(u)}(\Q(u))$.
 Since the points $P_1$ and $P_2$ have non-constant coordinates, it follows that the points are of infinite order. Since $\deg(g)=6$, one knows that the Mordell-Weil rank $r$ cannot be greater than $2$, see Remark \ref{rem1}.
\end{Proof}

The following lemma describes the Mordell-Weil group of some elliptic curve that will be used towards proving our main result.
\begin{Lemma}
\label{lem:auxiliary}
Let $\displaystyle \alpha\in\Q-\{0,\pm1\}$. The algebraic curve described by \[C_{\alpha}:y^2=\alpha^2x^4-(1+\alpha^2)^2x^2+4\alpha^2\] is an elliptic curve of positive Mordell-Weil rank over $\Q$.
\end{Lemma}
\begin{Proof}
Since $\displaystyle (x,y)=\left(0,2\alpha\right)\in C_{\alpha}(\Q)$, the algebraic curve $C_{\alpha}$ is an elliptic curve. Further, a Weierstrass equation describing $C_{\alpha}$, see \cite{Cremona}, is given by \[C'_{\alpha}:y^2=x^3-27(48\alpha^4+(1+\alpha^2)^4)x-54(1+\alpha^2)^2((1+\alpha^2)^4-144\alpha^4)\] with $(X,Y)\in C'_{\alpha}(\Q)$ where
\begin{eqnarray*}
X&=&\frac{3}{4\alpha^4}(1+\alpha^2)^2(3+12\alpha^2-22\alpha^4+12\alpha^6+3\alpha^8),\\
Y&=&\frac{27}{8\alpha^6}(-1+\alpha^2)^2(1+11\alpha^2+37\alpha^4+47\alpha^6+47\alpha^8+37\alpha^{10}+11\alpha^{12}+\alpha^{14}).
\end{eqnarray*}
A simple specialization argument yields that $(X,Y)$ is a point of infinite order in $C_{\alpha}'(\Q)$.
\end{Proof}

\section{Rank 2 Quadratic twists of pairs of Legendre curves}
\label{sec:5}

The following theorem contains our main result. An infinite family of pairs of Legendre curves is introduced together with a polynomial of degree $12$ such that the simultaneous quadratic twists of the pair by this polynomial are both of rank at least $2$. In fact, the quadratic twists are parametrized by the elliptic curve introduced in Lemma \ref{lem:auxiliary}.

\begin{Theorem}
\label{thm:family1}
Let $\alpha\in\Q-\{0,\pm1\}$. Consider the two Legendre elliptic curves
\begin{eqnarray*}
E_1&:& y^2=f_1(x)=x(x-1)\left(x+\frac{(\alpha^2+1)^2}{(\alpha^2-1)^2}\right)\\
E_2&:&y^2=f_2(x)=x(x-1)\left(x-\frac{(\alpha^2+1)^2}{4\alpha^2}\right).
\end{eqnarray*}
Set $g_{\alpha}(t)=(t^4+4)(t^4(\alpha^2-1)^2+4t^2(\alpha^2+1)^2
+4(\alpha^2-1)^2)(\alpha^2t^4 -(\alpha^2+1)^2t^2 +4\alpha^2)$, where $t$ is the $t$-coordinate of a rational point on the curve $u^2=\alpha^2t^4 -(\alpha^2+1)^2t^2 +4\alpha^2$. Then the quadratic twists $E_1^{g_{\alpha}(t)}$ and $E_2^{g_{\alpha}(t)}$ are of Mordell-Weil rank $r\ge 2$ over $\Q$. Moreover, the points $$\left(\frac{(t^4+4)(\alpha^2+1)^2}{4u^2},\frac{(t^2-2)(\alpha^2+1)^3}{8(\alpha^2-1)u^4}\right),
\left(\frac{t^2}{4}+\frac{1}{t^2},\frac{t^2-2}{8t^3(\alpha^2-1)u}\right)$$
are two independent points of infinite order in $E_1^{g_{\alpha}(t)}(\Q)$, and the points
 $$\left(\frac{\left(t^4+4\right)\left(\alpha ^2+1\right)^2 }{4 u^2},\frac{t\left(\alpha ^2+1\right)^3}{8\alpha u^4}\right),
\left(\frac{t^4+4}{\left(t^2+2\right)^2},\frac{t}{ (t^2+2)^3\alpha u}\right)$$
 are two independent points of infinite order in $E_2^{g_{\alpha}(t)}(\Q)$.
\end{Theorem}
\begin{Proof}
Checking that the given points are points on the corresponding elliptic curve is a direct calculation. That the points are of infinite order can be seen either by observing that they have non-constant coordinates or by specialization. That the points are independent can be checked by noticing that the automorphism $u\mapsto -u$ of $\Q(u)$ fixes the first point and sends the second point to its inverse.
\end{Proof}
In the following remark we illustrate the method that enabled us produce the elliptic curves and the rational points in Theorem \ref{thm:family1}.
\begin{Remark}
\label{rem:main}
 We set $\displaystyle \lambda_1=-\frac{(\alpha^2+1)^2}{(\alpha^2-1)^2}$ and $\displaystyle \lambda_2=\frac{\lambda_1}{1+\lambda_1}=\frac{(\alpha^2+1)^2}{4\alpha^2}$, $\alpha\in\Q-\{0,\pm1\}$. Moreover, we set $\displaystyle z=\frac{-\lambda_2+\lambda_1T^2}{-1+T^2}$ which implies that $f_2(z)=T^2 f_1(z)$.

In Theorem \ref{Thm:1new}, we take $\displaystyle h_i(z)=\frac{\lambda_i z}{(\lambda_i+1)z-\lambda_i}$, $i=1,2$. In fact, $h_i$ is the linear fractional transformation which sends the zeros $0,1,\lambda_i$ of $f_i$ to $0,\lambda_i,1$, respectively. The linear polynomial $k_i(z)=\lambda_i\left((1+\lambda_i)z-\lambda_i\right)$ is the square-free part of $f_i(h_i(z))/f_i(z)$.

Now, taking $\displaystyle T= \frac{(\alpha^2-1)t}{\alpha(t^2-2)}$, one can easily check that the genus one curve $C$ in Theorem \ref{Thm:1new} has a rational point if $\alpha^2t^4-(1+\alpha^2)^2t^2+4\alpha^2$ is a rational square. Therefore, setting $C_{\alpha}$ to be the elliptic curve $u^2=\alpha^2t^4-(1+\alpha^2)^2t^2+4\alpha^2$, one knows that the existence of a point $(t,u)$ in $C_{\alpha}(\Q)$ together with the fact that $f_2(z(t))/f_1(z(t))$ is the rational square $T^2$ yield the existence of a rational point on the genus one curve
\[C:z_1^2=k_1(z(t)),\;z_2^2=k_2(z(t)),\;z_3^2=f_2(z(t))/f_1(z(t)),\] namely,
$\displaystyle(z_1,z_2,z_3)=\left(\frac{\left(\alpha ^2+1\right)^3 t}{\left(\alpha ^2-1\right)^2 u},\frac{\left(\alpha ^2+1\right)^3 \left(t^2+2\right)}{8 \alpha ^2 u},\frac{\left(\alpha ^2-1\right) t}{\alpha  \left(t^2-2\right)}\right)$.
Since $C_{\alpha}$ is an elliptic curve of positive rank, Lemma \ref{lem:auxiliary}, then $C$ is also an elliptic curve of positive rank. Consequently, $C(\Q)$ contains infinitely many rational points. Now, according to Theorem \ref{Thm:1new}, one obtains that $\rank E_1^{f_1(z(t))}(\Q)\ge 2$ and $\rank E_2^{f_1(z(t))}(\Q)\ge 2$, where $t$ is the $t$-coordinate of a point in $C_{\alpha}(\Q)$ but for finitely many exceptions. In particular, Theorem \ref{thm:family1} provides us with infinitely many rational numbers $d$ such that the quadratic twists $E_1^d$ and $E_2^d$ are of rank at least $2$ over $\Q$.
\end{Remark}

One can repeat the procedure in Remark \ref{rem:main} choosing different linear fractional transformations in order to obtain a new family of pairs of elliptic curves $E_1$ and $E_2$ over $\Q$ and infinitely many rationals $d$ such that $E_1^d$ and $E_2^d$ are of Mordell-Weil rank at least $2$ over $\Q$. In fact, by taking $\displaystyle h_{i}=\frac{t-\lambda_i}{(2-\lambda_i)t-1}$, $i=1,2$, in Remark \ref{rem:main}, one may obtain the following result.

\begin{Theorem}
Let $\alpha\in\Q-\{0,\pm1\}$. Consider the two Legendre elliptic curves
\begin{eqnarray*}
E_1&:& y^2=f_1(x)=x(x-1)\left(x-\frac{2(1+\alpha^4)}{(-1+\alpha^2)^2}\right)\\
E_2&:&y^2=f_2(x)=x(x-1)\left(x+\frac{(-1+\alpha^2)^2}{4\alpha^2}\right).
\end{eqnarray*}
Set \[g_{\alpha}(t)=-(t^4+4)((\alpha^2-1)^2t^4+4(\alpha^2+1)^2t^2
+4(\alpha^2-1)^2)(\alpha^2t^4 -(\alpha^2+1)^2t^2 +4\alpha^2),\] where $t$ is the $t$-coordinate of a rational point on the curve $u^2=\alpha^2t^4 -(\alpha^2+1)^2t^2 +4\alpha^2$. Then the quadratic twists $E_1^{g_{\alpha}(t)}$ and $E_2^{g_{\alpha}(t)}$ are of Mordell-Weil rank $r\ge 2$ over $\Q$. Moreover, the points $$\left(-\frac{t^4(-1+\alpha^2)^2+4t^2(1+\alpha^2)^2+4(-1+\alpha^2)^2}{4u^2},
\frac{(t^2-2)(\alpha^2+1)^3}{8(\alpha^2-1)u^4}\right),
 \left(-\frac{(t^2-2)^2}{4 t^2},\frac{t^2-2}{8t^3(\alpha^2-1)u}\right)$$
are two independent points of infinite order in $E_1^{g_{\alpha}(t)}(\Q)$, and the points
 $$\left(-\frac{t^4(-1+\alpha^2)^2+4t^2(1+\alpha^2)^2+4(-1+\alpha^2)^2}{4u^2},
\frac{t(\alpha^2+1)^3}{8\alpha u^4}\right),
\left(\frac{4 t^2}{(t^2+2)^2},\frac{t}{(t^2+2)^3\alpha u}\right)$$
 are two independent points of infinite order in $E_2^{g_{\alpha}(t)}(\Q)$.
\end{Theorem}

It is worth mentioning that we can obtain further families by changing our choices of $h_1$ and $h_2$, yet the resulting polynomial $g_{\alpha}$ has large coefficients.

\hskip-18pt\emph{\bf{Acknowledgements.}}
We would like to thank the referee for many comments, corrections, and
suggestions that helped the authors improve the manuscript significantly.

\end{document}